\documentclass[11pt,a4paper]{article}
\usepackage[latin1]{inputenc}
\usepackage[american]{babel}      
\usepackage[T1]{fontenc}
\usepackage{amsmath,amssymb,amsthm}
\usepackage{graphicx}
\usepackage{longtable}
\usepackage[export]{adjustbox}

\title{Detection of thermal bridges\\ from thermographic images for the analysis\\ of buildings energy performance}
         
\author{{\bf Francesco Asdrubali}$^1$, \hskip0.05cm {\bf Giorgio Baldinelli}$^2$, \hskip0.05cm {\bf Francesco Bianchi}$^2$, \\ {\bf Danilo Costarelli}$^3$,\hskip0.2cm {\bf Antonella Rotili}$^2$,\\ \hskip0.4cm {\bf Marco Seracini}$^3$ \hskip0.3cm and \hskip0.3cm {\bf Gianluca Vinti}$^3$  \thanks{{\footnotesize 1: Department of Engineering, University of Roma Tre, Via V. Volterra 62, 00146, Roma, Italy, {\tt 	francesco.asdrubali@uniroma3.it}; 2: Department of Engineering, University of Perugia, Via G. Duranti 67, 06125, Perugia, Italy, {\tt giorgio.baldinelli@unipg.it}, {\tt bianchi.unipg@ciriaf.it}, {\tt rotili.unipg@ciriaf.it}; 3: Department of Mathematics and Computer Science, University of Perugia, Via Vanvitelli 1, 06123 Perugia, Italy,
 {\small {\tt danilo.costarelli@unipg.it}, {\tt marco.seracini@unipg.it}}, {\small {\tt gianluca.vinti@unipg.it}} } } }

\date{}

\newcommand{\miu}{\leq}

\newcommand{\N}{\mathbb{N}}
\newcommand{\Z}{\mathbb{Z}}
\newcommand{\R}{\mathbb{R}}
\newcommand{\uu}{\underline{u}}

\newcommand{\xx}{\underline{x}}

\newcommand{\kk}{\underline{k}}

\newcommand{\be}{\begin{equation}}
\newcommand{\ee}{\end{equation}}

\newtheorem{definition}{Definition}[section]

\newtheorem{theorem}[definition]{Theorem}

\begin{document}

\maketitle  

\begin{abstract}
In this paper, we develop a procedure for the detection of the contours of thermal bridges from thermographic images, in order to study the energetic performance of buildings. Two main steps of the above method are: the enhancement of the thermographic images by an optimized version of the mathematical algorithm for digital image processing based on the theory of sampling Kantorovich operators, and the application of a suitable thresholding based on the analysis of the histogram of the enhanced thermographic images. Finally, an accuracy improvement of the parameter that defines the thermal bridge is obtained.
\vskip0.3cm
\noindent
  {\footnotesize AMS 2010 Classification: 65D15, 65T99, 65C99, 41A30}
\vskip0.1cm
\noindent
  {\footnotesize Key words and phrases: Sampling Kantorovich operators, approximation results, thermographic images, image processing.} 
\end{abstract}

\section{Introduction} \label{sec1} 

The thermographic survey on the building envelope is a useful non-invasive method to detect thermal bridges that reduce the overall energy performance of buildings. 

In \cite{ASBA1,ASBA2} the above investigation has been performed by the help of a suitable index, the so-called {\em incidence factor of thermal bridge $I_{tb}$}, which points out the energetic incidence of a thermal bridge on the basis of the temperature decrement that it causes. The accuracy of this analysis depends on various aspects. One of the most important relies in the correct detection of the areas belonging to thermal bridges; such analysis described in \cite{ASBA1,ASBA2} is strictly operator dependent. A further important aspect is given by the quality (i.e., the resolution) of the thermographic data available to perform the mentioned energetic analysis.

    In the present paper, we introduce a segmentation procedure which allows to detect the thermal bridges from thermographic images, becoming the energy analysis of the buildings automatic and more accurate than the original one proposed in \cite{ASBA1}. 
    
    The algorithm developed in this paper is characterized by various steps, based on mathematical method of approximation theory and on techniques of Digital Image Processing (D.I.P.). Firstly, the thermographic images are reconstructed and enhanced in their resolution, by the application of the sampling Kantorovich (S-K) algorithm, see e.g., \cite{COVI2,CLCOMIVI1,CLCOMIVI2}. The latter method can be deduced from some approximation results concerning the theory of the well-known S-K operators, which has been deeply investigated in last years, see e.g., \cite{COVI1,COVI4,COVI5,COMIVI1,ORTA1}. In particular, here we developed a suitable numerical optimized version of the S-K algorithm, useful to solve the application problem treated in this paper. Note that the S-K algorithm has been firstly introduced in \cite{COVI2}, but its original implementation was not-optimized and required a longer CPU-time of execution.  
    
    The sampling Kantorovich operators represent in practice an $L^1$-version of the classical generalized sampling operators introduced by P.L. Butzer and his school at Aachen in the `80s, with the aim to obtain an approximate sampling formula. In some recent papers, it has been proved that the S-K operators are suitable for applications to real-world cases based on the study of suitable images; for instance, some models for the study of the behavior of buildings under seismic action have been successfully obtained by the applications of the S-K algorithm for the enhancement of thermographic images, see e.g., \cite{CLCOMIVI1,CLCOMIVI2}.
    
    Here, the S-K algorithm for image enhancement has been implemented by a bivariate Jackson-type kernel, and the original resolution of thermographic images has been improved.
    
    Then, by a probabilistic interpretation of the histogram of the pixels associated to the above enhanced thermographic images (see e.g., \cite{KASAWO1}), we determined a suitable threshold value which can be used in order to segment the thermal bridge.
    
    The validation of the proposed method has been obtained from the experimental results effected in an ad-hoc built hot-box setup with controlled laboratory conditions (see e.g., \cite{ABBP1,BBPR1}). More precisely, two types of bi-dimensional thermal bridges, with different shape, have been built and successively tested. The numerical results show that the proposed algorithm, other to identify the geometry of the thermal bridge generated in walls composite by different materials, allows to improve the energy analysis of the buildings with respect to the original approach given in \cite{ASBA1}.
    
   Indeed, such improvement has been validated by a comparison among the factor $I_{tb}$, firstly computed according with the data detected by the probes in the hot-box, the original procedure developed in \cite{ASBA1}, and the method here proposed, i.e., working with the thermographic image, enhanced by the S-K algorithm with the shape of the thermal bridge extracted by the automatic thresholding procedure.   
    
The above numerical results show that, in thermal bridges caused by different materials, the method here introduced provides results closer to the most accurate approach, i.e., to those ones computed with the help of the probes, as numerically showed in Section \ref{sec5}.


\section{Approximation by sampling Kantorovich operators and applications to image processing} \label{sec2}

In this section, we give a background concerning the main theoretical and applications aspects of the theory of sampling Kantorovich (S-K) operators, see e.g., \cite{BABUSTVI2,VIZA1} for results concerning functions of one variable, and \cite{COVI1,COVI4,COVI5} for what concerns functions of several variables. 

  The above family of operators are typically used in approximation theory in order to reconstruct not-necessarily continuous signals, such as images, see \cite{BURIST2,GOWO,COVI2}, and they are useful for image reconstruction and enhancement, see \cite{CLCOMIVI1,CLCOMIVI2}. 
  
   We begin, by recalling the definition of the kernel functions used in order to define the above approximation process. 
  
  In what follows, we will define as kernel any multivariate function $\chi: \R^n \to \R$, which satisfies the following conditions:
\begin{itemize}
\item[($\chi 1$)] $\chi$ is summable on $\R^n$, and bounded in a ball containing the origin of $\R^n$;

\item[($\chi 2$)] For every $\xx \in \R^n$:
$$ 
\sum_{\kk \in \Z^n}\chi(w \xx- \kk) = 1;
$$

\item[($\chi 3$)] For some $\beta>0$, we assume that the discrete absolute moment of order $\beta$ is finite, i.e.,
$$
	m_{\beta}(\chi)\ :=\ \sup_{\uu \in \R}\, \sum_{\kk \in \Z^n}\left|\chi(\uu-\kk)\right|\cdot\|\uu-\kk\|^{\beta}\ <\ +\infty.
$$
\end{itemize}

\begin{figure}
\centering
\includegraphics[scale=0.433]{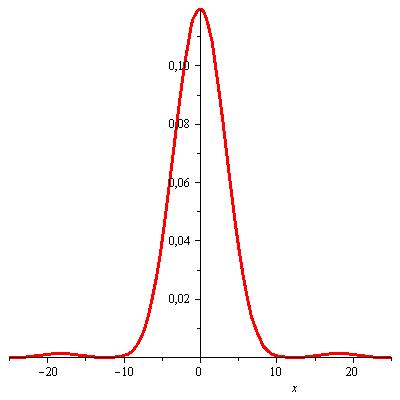}
\includegraphics[scale=0.433]{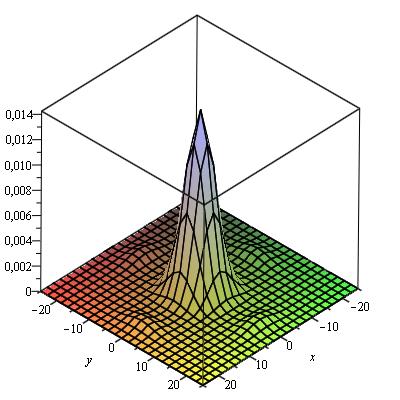}
\caption{{\small The plots of the one-dimensional and the bivariate Jackson-type kernels $J_2(x)$ (left) and ${\cal J}^2_{2}(\xx)$ (right), with $\alpha=1$.}} \label{JAK}
\end{figure}

\begin{figure}
\centering
\includegraphics[scale=0.433]{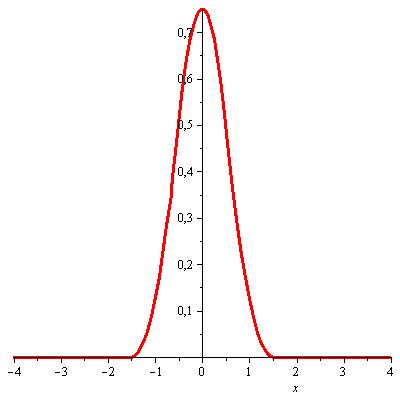}
\includegraphics[scale=0.433]{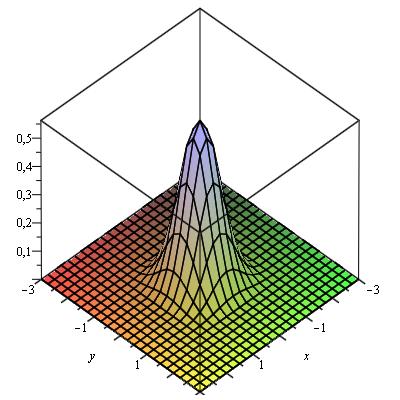}
\caption{{\small The plots of the one-dimensional and the bivariate central B-spline $M_3(x)$ (left) and ${\cal M}^2_3(\xx)$ (right), respectively.}} \label{m3}
\end{figure}

  We immediately provide some typical examples of kernels which satisfy all the above assumptions $(\chi 1)$, $(\chi 2)$, and $(\chi 3)$. The most used method to construct multivariate kernels is to consider the product of $n$ kernels of one variable, see \cite{BUFIST,COGA1,COGA2,CO2}. Indeed, for instance, the definition of the multivariate Fej\'{e}r kernel can be formulated as follows:
\be
\mathcal{F}_n(\xx)= \prod^n_{i=1}F(x_i), \hskip1cm \xx = (x_1, ..., x_n) \in \R^n,
\ee
where $F(x)$, $x \in \R$, denotes the univariate Fej\'{e}r kernel, which is defined by:
\be
F(x)\ :=\ \frac{1}{2}\, \mbox{sinc}^2\left(\frac{x}{2}\right), \hskip1cm x \in \R,
\ee
where the well-known $sinc$-function is that defined as $\sin (\pi x)/\pi x$, if $x \neq 0$, and $1$ if $x=0$, see e.g., \cite{ANVI1,BABUSTVI3,TA1,VIZA3,AGBA}. By the $sinc$-function it is possible to define another class of kernel, which is widely used, i.e., the Jackson-type kernels, see e.g., \cite{BUNE,BAMUVI,IC,CLCOMIVI1,CLCOMIVI2}.
The multivariate expression of the Jackson-type kernels (see e.g., Fig. \ref{JAK}), is the following: 
\be
{\cal J}^n_k(\xx)\ :=\ \prod^n_{i=1}J_k(x_i), \hskip1cm \xx=(x_1,...,x_n) \in \R^n,
\ee
where $J_k(x)$, $x\in \R$ are defined by:
\be
J_k(x)\ :=\ c_k\, \mbox{sinc}^{2k}\left(\frac{x}{2k\pi\alpha}\right), \hskip1cm x \in  \R,
\ee
with $k \in \N$, $\alpha \geq 1$, and $c_k$ is a non-zero normalization coefficient, given by:
$$
c_k\ :=\ \left[ \int_{\R} \mbox{sinc}^{2k}\left(\frac{u}{2 k \pi \alpha} \right) \, du \right]^{-1}.
$$
Since the $sinc$-function has unbounded support, we usually say that $F(\xx)$ and ${\cal J}^n_k(\xx)$ are not duration limited kernels. To know the duration of a kernel is important in order to implement the numerical evaluation of the operators studied in this section. In fact, operators based upon kernels with unbounded duration, need to be truncated for the evaluation. For the latter reason, we also provide examples of duration limited kernels. For instance, we can consider the well-know central B-spline of order $s$ (see e.g., \cite{BUNE,CHDI1,LIWA1,KV1,SA1}), defined by:
\be
M_s(x)\ :=\ \frac{1}{(s-1)!}\, \sum^s_{i=0}(-1)^i \left(\begin{array}{l} \!\! 
s\\
\hskip-0.1cm i
\end{array} \!\! \right)
\left(\frac{s}{2}+x-i\right)^{s-1}_+,
\ee
where the function $(x)_+ := \max\left\{x,0\right\}$ denotes the positive part of $x \in \R$. The corresponding multivariate spline kernels (see e.g., Fig. \ref{m3}) are then defined by:
\be
\mathcal{M}^n_s(\xx)\ :=\ \prod^n_{i=1}M_s(x_i), \hskip1cm \underline{x}=(x_1,...,x_n) \in \R^n.
\ee 
For others useful examples of kernel functions, see e.g., \cite{BUFIST,COVI6,COVI7,COVI8}.

  Now, we are able to recall the definition of the multivariate sampling Kantorovich operators (\cite{COVI1}), of the form:
\be \label{KANTOROVICH}
(S_w f)(\xx)\ :=\ \sum_{\kk \in \Z^n} \chi(w\xx-\kk)\, \left[ w^n \int_{R_{\kk}^w }f(\uu)\ d\uu \right], \hskip0.7cm \xx \in \R^n,
\ee
where $w>0$, and $f: \R^n \to \R$ is a locally integrable function such that the above series is convergent for every $\xx \in \R^n$, and 
$$
R_{\underline{k}}^w\ :=\ \left[\frac{k_1}{w},\frac{k_1+1}{w}\right]\times\left[\frac{k_2}{w},\frac{k_2+1}{w}\right] \times ... \times \left[\frac{ k_n }{w},\frac{ k_n+1 }{w}\right],
$$
are the sets in which we consider the mean values of the signal $f$.

For the sake of completeness, we must observe that the above operators $S_w$ can be reviewed as an $L^1$-version of the generalized sampling operators, introduced by P.L. Butzer and his school in 1980s, and widely studied in last thirty years, see e.g., \cite{BUFIST,BURIST2,KITA1,KIME1}. The sampling Kantorovich operators, as well as the generalized sampling operators, are a kind of quasi-interpolation operators, see e.g., \cite{CHDI1,VE1,TA1,SK}.

  We now recall the main approximation results for the S-K operators generated by the kernel $\chi$ satisfying $(\chi i)$, $i=1,2,3$, on which the algorithm for image reconstruction and enhancement is based.
\begin{theorem}[\cite{COVI1}] \label{theorem1}
Let $f: \R^n \to \R$ be a given bounded signal. Then:
$$
\lim_{w \to +\infty} (S_wf)(\xx)\ =\ f(\xx),
$$
at any point of continuity of $f$. Moreover, if $f$ is uniformly continuous on $\R^n$, it turns out that:
$$
\lim_{w \to +\infty} \| S_wf - f\|_{\infty}\ =\ \lim_{w \to +\infty} \sup_{\xx \in \R^n}  | (S_wf)(\xx) - f(\xx) |\ =\ 0.
$$
Finally, if the signal $f$ belongs to $L^p(\R^n)$, $1 \miu p < +\infty$, we have:
$$
\lim_{w \to +\infty}\| S_wf - f\|_p\ =\ \lim_{w \to +\infty} \left( \int_{\R^n} | (S_wf)(\xx) - f(\xx) |^p\, d\xx  \right)^{1/p}\ =\ 0.
$$
\end{theorem}
Note that, the proof of Theorem \ref{theorem1} for the multivariate S-K operators given in \cite{COVI1} are constructive. 

   For what concerns the result in the $L^p$-setting, this can be deduced by a general theorem proved in Orlicz spaces, see e.g., \cite{MU1,BAMUVI,COVI9}.

  Now, we point out the main steps needed to implement the algorithm for image reconstruction and enhancement based on the sampling Kantorovich operators. 

   First of all, we recall that a gray scale image $A$ (matrix) of size $n \times m$ can be modeled as a step function $I$, belonging to $L^p(\R^2)$, $1 \leq p <+\infty$, as follows:
\be
I(x,y)\ :=\ \sum^{n}_{i=1}\sum^{m}_{j=1}a_{ij} \cdot \textbf{1}_{ij}(x,y), \hskip1cm (x,y) \in \R^2,
\ee
where ${\bf 1}_{ij}(x,y)$, $i,j =1,2,...,m$, are the characteristics functions of the sets $(i-1,\ i]\times(j-1,\ j]$, see e.g., \cite{ZH1,CHYUVE1,GOWO,KI}. 
Then, as stated by Theorem \ref{theorem1}, the family of the S-K operators $(S_w I)_{w>0}$ (for some kernel $\chi$ in the two-dimensional case) approximates $I$, pointwise at any points belonging to the interior of the sets $(i-1,\ i]\times(j-1,\ j]$, $i,\, j = 1, ..., m$, and in $L^p$-sense. 

  Now, in order to obtain a new image (matrix) that approximates the original one, it is sufficient to sample $S_w I$, for some $w>0$, with a fixed sampling rate. In particular, we can reconstruct the approximating images taking into consideration different sampling rates. If we choose an higher sampling rate, as a result of the approximation process, we get a new image that has a better resolution than the original one. Therefore, the above procedure is also useful for image enhancement. The algorithm has been implemented by using MATLAB, and a pseudo-code has been reported in Table \ref{table1}. Examples of reconstructed image by the above procedure can be found, e.g., in \cite{CLCOMIVI1}.
\begin{table}[!h]
\centering
\begin{tabular}{c}
\hline 
\\ 
\underline{Objective}:  Reconstructing and improving\\
\hskip1.5cm the resolution of the original \\
\hskip1.8cm bivariate image $I$ by sampling\\ 
\hskip1.6cm Kantorovich operators based\\
\hskip1.4cm upon the bivariate kernel $\chi$. \\
\\
\hskip0.3cm \underline{Inputs}: Original image $I$ ($n \times m$ pixel\\
\hskip0.9cm resolution), the parameter\\ 
\hskip1.8cm$w>0$ and the scaling factor $R$.
\\
\\
\hskip0.5cm $\bullet$ Choice and definition of the kernel\\ \hskip-2.8cm function $\chi$;
    \\
\hskip0.1cm $\bullet$ Size of the reconstructed image:\\ 
\hskip-2.0cm ($n \cdot R$)$\times$($m \cdot R$);
\\
$\bullet$ Computation of matrices of the\\ 
\hskip0.3cm mean values (samples) of $I$ by\\ 
\hskip0.4cm means of the Kronecker matrix\\ 
\hskip-3.2cm product.
\\
$\bullet$ Definition of the vectors conta-\\
\hskip-0.4cm ining the arguments of $\chi$.
\\
\\
\hskip0.0cm \underline{Iteration}: Summation over $\kk$ of all non\\ 
\hskip0.6cm zero terms of the form\\ 
\hskip1.8cm $\chi(w \underline{x} - \kk) \cdot \left[ w^2 \int_{R_{\underline{k}}^w} I(\underline{u})\, d\underline{u} \right]$, \\
\hskip1.2cm for a
suitable fixed grid of\\
\hskip-1.6cm points $\underline{x}$.
\\
\\
\hskip0.0cm \underline{Output}: The reconstructed image of\\ 
\hskip1.3cm resolution
($n \cdot R$)$\times$($m \cdot R$).
\\
\\
\hline 
\end{tabular}
\caption{Pseudo-code of the sampling Kantorovich algorithm for image reconstruction and enhancement.}
\label{table1}
\end{table}
%


\section{The thermographic data and the segmentation procedure} \label{sec3} 

In this section we present the data used for the analysis 
of the thermal bridges. 
In civil and energy engineering, in order to effect non-invasive investigations, infrared thermography images are usually used. 
The thermography is a technique which allows to measure the heat flux associated with the infrared radiation emitted from every object without direct contact, then it supplies a non-invasive technique for investigating buildings, see e.g., \cite{CACA1,BOCAME1,BOCAME2}.
The thermography exploits the peculiarity that all the objects (having a certain temperature higher than the absolute zero) emit radiation in the infrared range (wave length of $700\,\mathrm{nm}\textrm{--} 1\,\mathrm{mm}$, which corresponds to frequencies of $430\,\mathrm{THz}\textrm{--} 300\,\mathrm{GHz}$), which is located between the visible radiation (in particular the red component) and the microwave range.

The result of a thermographic survey is a bi-dimensional thermal mapping of the heat flux of the object expressed in temperature when the emissivity is known.

  Thermal bridges are zones of buildings that present a thermal flow higher than the adjacent constructive elements. For this reason, thermographic images appear to be the most appropriate tool to study such a problem, and to effect the analysis of the energetic performance of the building envelope.  
   
  The presence of the zones described above determines the energy losses of efficiency in buildings and it brings impacts on structural and comfort aspects.

Because of the higher heat flux across the walls, the areas affected by thermal bridges appear colder than those nearest, i.e., darker in the fictitious visual representation provided by the thermal camera when an inner wall is analyzed in winter conditions.

  We studied two cases of thermal bridges, depending on the morphology of the walls on which they appear:
\begin{enumerate}
\item 2D thermal bridges - pillar;
\item 2D thermal bridges - beam-pillar joint;
\end{enumerate}
The first type of thermal bridges occurs on plain walls and develops in a straight direction, while the second type occurs again on plain walls but draws an angle on them.

\begin{figure}[!htb]
\hskip-0.9cm
  \includegraphics[scale=0.37]{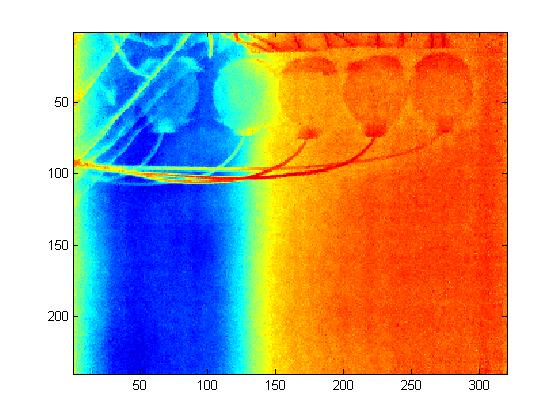}
  \hskip-0.7cm
  \includegraphics[scale=0.37]{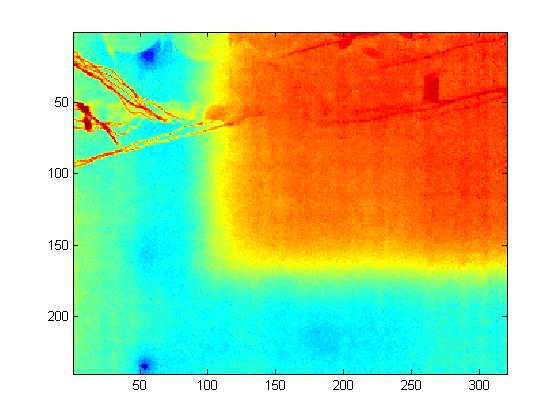}
  \caption{The thermal bridges built in a hot-box: on the left, a pillar 2D thermal bridge; on the right, a beam-pillar-joint 2D thermal bridge.}\label{fig1}
\end{figure}

In order to effect a complete analysis on thermal bridges and their energy performance, and to make the present dissertation as complete as possible, we built in hot-box setup concrete examples of each one of the thermal bridges above described, see Fig. \ref{fig1}. Having constructed in a laboratory such thermal bridges allows to know all their thermal characteristics. The latter fact will be used in order to evaluate the goodness of the analysis performed with the methodology here developed. 

   All data showed in Fig. \ref{fig1}, have been acquired by a Flir b360 thermal camera with a 320x240 pixel resolution.
The values retrieved from the thermal camera represent temperature data but they can be also interpreted as gray levels of an image, and they can be displayed in a readable format.

 In this paper the thermal values have been represented by images coded with a certain color-map instead of a gray-scale one.
The main scope of the present data analysis is to obtain an automatic procedure to extract the contours of the thermal bridge. Currently, it seems that no automatic procedures are available for the above task.

  At the aim of achieving this result, a segmentation procedure has been implemented and applied to the original matrix dataset $I=[N \times M]$, with N=320, M=240.
  To improve the accuracy of the border extraction, the sampling Kantorovich algorithm has been used on the native thermal data $I$, with a scaling factor  $R=2$, $w=15$, and by bivariate Jackson-type kernel generated by $J_{12}(x)$.  
In order to speed up the procedure, we approximated by truncation the S-K operators by neglecting all terms with values less than $10^{-4}$. This value has been computed taking into account the accuracy of the measurements (that here is $10^{-2}$), and other parameters. For a deeper discussion, see Section \ref{sec4}.

\begin{figure}[!htb]
\hskip-0.8cm
  \includegraphics[scale=0.37]{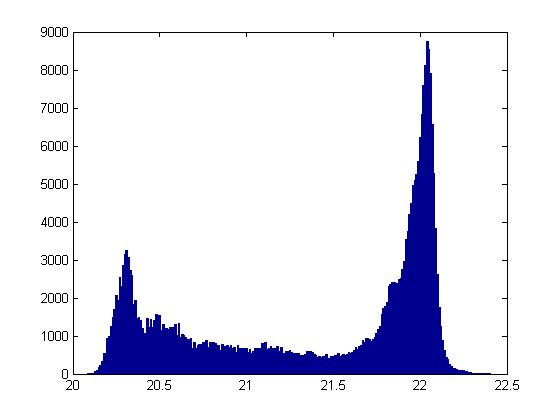}
  \hskip-0.6cm
  \includegraphics[scale=0.37]{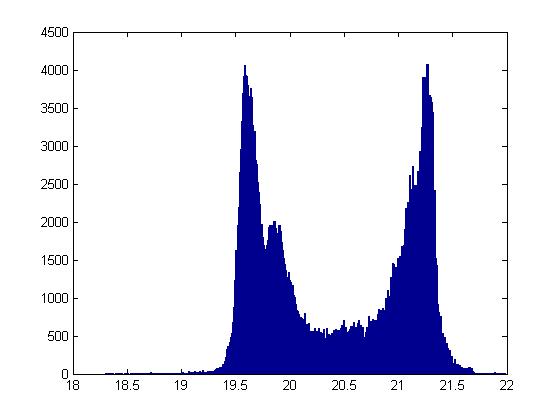}  \caption{Histograms corresponding to the thermal bridges given in Fig. \ref{fig1}, respectively, which have been processed by the S-K algorithm.}\label{fig5}
\end{figure}

  Once the data have been reconstructed using the S-K algorithm, the contours of the thermal bridges are extracted by a thresholding, using a parameter obtained from the associated histogram of temperatures (see Fig. \ref{fig5}). 

  The procedure for the identification of the threshold value is described in what follows. The assumption that the thermal bridges boundaries are defined as the zones where the temperature shows a significant gradient compared to its mean value on the entire surrounding area, has been made. 
  
\begin{figure}
\centering
\includegraphics[scale=0.62]{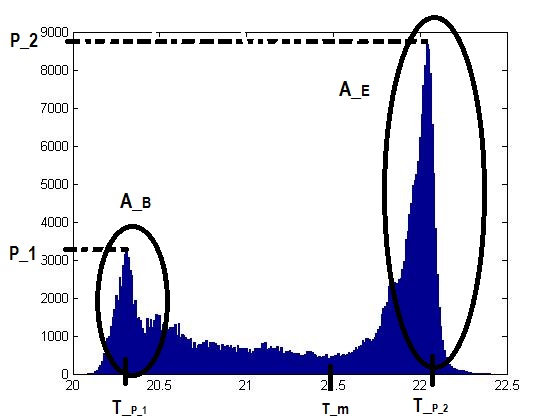}
\caption{{\small Identification of the various parameters on the histogram of the temperatures (concerning the first thermal bridge of Fig. \ref{fig1}), needed in order to compute the threshold value for the segmentation. The oval shaped curves have been inserted in the picture only to approximatively indicate those pixels belonging to the two bigger homogeneous areas, in fact denoted by $A_B$ and $A_E$.}} \label{fig3new}
\end{figure}

\begin{figure}[!htb]
\centering
\includegraphics[scale=0.45]{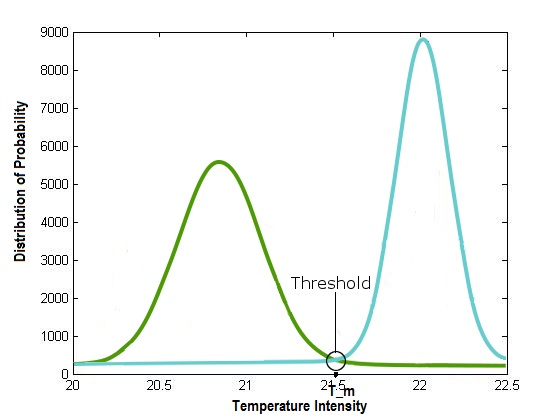}
\caption{Probabilistic interpretation of data values.}\label{fig6}
\end{figure}

  More precisely, by the above assumption, it is possible to verify as the data associated histogram exhibits two peaks corresponding to the two biggest relative maximum $P_1$ and $P_2$ ($P_1<P2$) respectively at coordinates $T_{P,1}$ and $T_{P,2}$, with an absolute minimum $T_{m}$ assumed on the interval $[T_{P,1},\, T_{P,2}]$, see e.g. Fig. \ref{fig3new}. We have chosen the threshold value $T_m$ as the temperature corresponding to the minimum between $P_1$ and $P_2$. This choice is made in view of the following probabilistic interpretation (see Fig. \ref{fig6}); indeed, $T_m$ can be associated to the minimum error due to the wrong classification of pixels located inside the thermal wall but classified as external, and viceversa.

  Note that, if the minimum is achieved at more than one value of temperature, we choose that one closer to the temperature corresponding to the higher value among $P_1$ and $P_2$.
  
  In Fig. \ref{fig3new}, $T_{P,1}$ and $T_{P,2}$ represent the values around which, points with homogeneous temperature are distributed; therefore, they highlight the thermal bridge area $A_B$ and the external homogeneous area $A_E$.
  
\begin{figure}
\hskip-0.6cm
\includegraphics[scale=0.35]{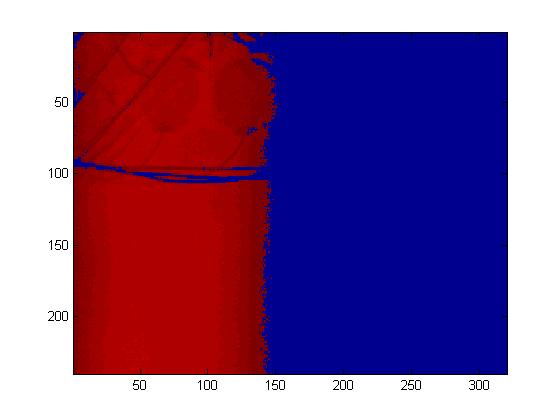}
\hskip-0.5cm
\includegraphics[scale=0.35]{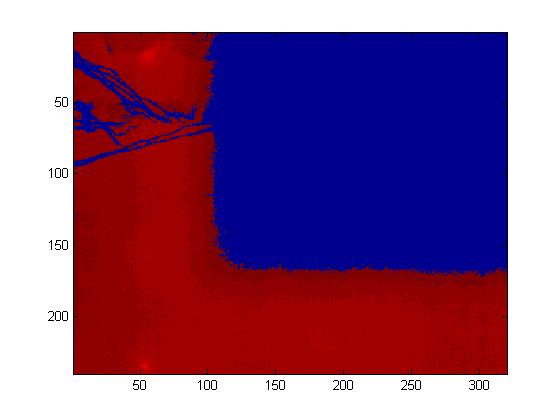}
\caption{Thermal bridges extracted using the threshold temperatures calculated by the algorithm. The homogeneous dark areas represent the zones outside the thermal bridges.}\label{fig7}
\end{figure}

  From the previous considerations, assuming  $P_1$ and $P_2$ to be the maximum of two distinct bell shaped functions which represent the distribution of probability of a point to belong to $A_B$ or $A_E$, is now possible to segment the original image $I$ choosing the threshold value equal to $T_{m}$. The results obtained from the above segmentation procedure concerning the thermal bridges of Fig. \ref{fig1}, have been shown in Fig. \ref{fig7}, and a pseudo-code of the algorithm for the determination of the threshold values $T_m$ is given in Table \ref{table11}.
\begin{table}[!h]
\centering
\begin{tabular}{c}
\hline 
\\ 
\underline{Objective}:  Segmentation of thermal bridges from\\
\hskip-1cm thermographic data\\ 
\\
\hskip0.8cm \underline{Inputs}: thermographic data of thermal bridges\\
\hskip-0.4cm of size $320\times 240$ points
\\
\\
\hskip2.5cm $\bullet$ Application of the sampling Kantorovich algorithm\\ 
\hskip0.4cm with scaling factor $R=2$ and $w=15$\\
\hskip0.2cm obtaining a matrix of size $640 \times 480$;
\\
\hskip2.5cm $\bullet$ Generation of the histogram associated to the data;
    \\
\hskip2.8cm $\bullet$ Computation of the minimum between $T_{P,1}$ and $T_{P,2}$.
\\
\\
\hskip-1.0cm \underline{Output}: The scalar threshold value $T_m$. 
\\
\\
\hline 
\end{tabular}
\caption{Pseudo-code of the segmentation algorithm for thermal bridge individuation.}
\label{table11}
\end{table}
\begin{table}[!h]
\begin{center}
\begin{tabular}{cccc}
\hline
\textbf{Thermal bridge number} & $T_{m}$ \\
\hline
TB. n. 1 & 21.50 $^{\circ}$C \\
TB. n. 2 & 20.36 $^{\circ}$C \\
\hline 
\end{tabular}
\caption{{\small Threshold values for the analyzed thermal bridges.}} \label{tab_threshold}
\end{center}                          
\end{table}
With the $T_{m}$ value is now possible to divide the data in two different zones, one individuating $A_B$, the other $A_E$. This method has been standardly used in literature for image segmentation (\cite{GOWO, KLME}) as well as for bimodal distribution statistical analysis. In Table \ref{tab_threshold}, the temperature of thresholding $T_{m}$, determined by using the described procedure, is given for each one of the two thermal bridges under analysis.
   
     The results obtained from the applications of the procedure described in this section are given in Fig. \ref{fig7}, where the shapes of the thermal bridges of Fig. \ref{fig1} have been shown. The energy behavior of the thermal bridges built in hot-box setup will be performed in Section \ref{sec5}.


\section{Numerical optimization of the S-K algorithm} \label{sec4}

\begin{figure}
\includegraphics[scale=0.32]{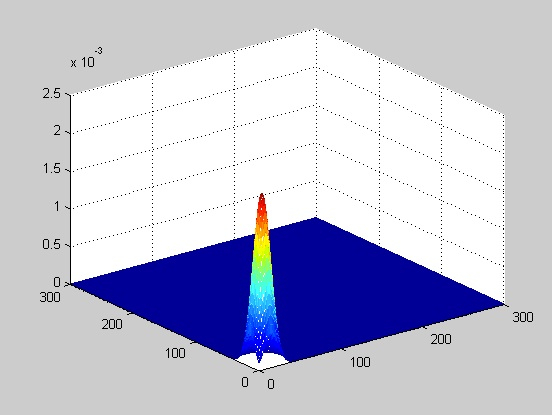}
\hskip0.1cm
\includegraphics[scale=0.32]{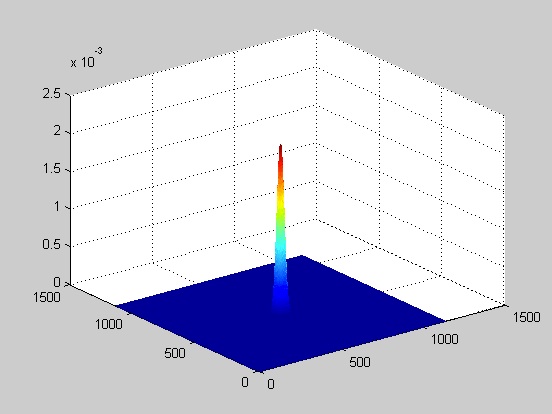}
\caption{Example of the bivariate Jackson's kernel ${\cal J}^2_{2}$: (left) calculated in a single step of the procedure (1), (right) computed as in the procedure (2).}\label{fig8}
\end{figure}

As described in Section \ref{sec2}, the reconstruction procedure, when implemented according to Table \ref{table1}, essentially reduces itself to multiply the chosen 2-dimensional kernel, computed in a suitable grid of nodes, for the starting data set, i.e., the gray-levels image $I$.
This numerical calculation can be operated following two basically different approaches:
\begin{itemize}
\item[(1)] an implementation of the previous theory, where in order to compute each value of the reconstructed data set, we recalculated the kernel matrix (see Fig. \ref{fig8}, left);

\item[(2)] first, we calculate the kernel matrix, sized according to $w$, at all the nodes necessary to complete the whole procedure, and then we select at each step the portion of the matrix of interest (see Fig. \ref{fig8}, right). Moreover, we adopted a numerical truncation of the above kernel matrix, in oder to neglect all the values that, certainly, are not significant in the final computation of the reconstructed image. This is possible, since the main part of the kernel (that provides a concrete contribution to the reconstruction procedure) is concentrated in a small portion of the domain (see Fig. \ref{fig8} again). A detailed description of the above optimization procedure follows below.
\end{itemize}
Each of these approaches has its advantage: the item (1) requires less memory occupation than the item (2), while item (2) requires a significant small number of arithmetic operations than (1).

    Indeed, in case (1), given a gray-scale image of size $N \times M$, the kernel memory occupation can be estimated by:
$$
N\, M\, w^2\, B,
$$
where $B$ is the number of bits used for the representation of every data value. In item (2) we have instead:
$$
N\, M\, w^2\, R^2\, B,
$$
where $R$ is the scaling factor. It is clear that if $R>1$, for the same value of $w$, $N$, $M$ and $B$, the approach (2) is more expensive in therms of memory occupation than (1).

   On the other hand, since procedure (2) computes the kernel matrix only once, and moreover it is also characterized by a numerical optimization (truncation), it requires a small number of arithmetic operations than the procedure (1).
   
   More precisely, the numerical truncation given in item (2) consists of determining a threshold value $\overline{k}>0$, useful to identify all the elements of the kernel matrix that can be neglected, since they are not significant in the computation of the output.   
   
   Recalling that, the previous thermographic data have been recorded by a thermal camera whose measurament resolution is $P:=10^{-2}$, the threshold value $\overline{k}$ can be obtained by the following:
$$
\overline{k}\ =\ 	\frac{4\, 10^{-1}\, P}{w^2\, N\, M\, A},
$$
where:
$$
A\ :=\ \max\left\{a_{i,j}:\ i=1, ..., N,\ j=1, ..., M \right\}, \hskip1cm I=(a_{i,j})_{i, j}.
$$
In this way, we neglected all values which contribute (in the worse case) of an amount less or equal to $0.4\, P$.

  With the purpose of establishing the above proposed procedures is better for our aim, we estimated the CPU-time for the reconstruction of some data sets on a computer system equipped with i7 quad core CPU and 8 GB of ram, running Matlab\textcircled c ver. R2014b on Windows\textcircled c  OP Windows 7\textcircled c, Service pack 1, 64-bit parallelism. We have considered the rate of growth of the CPU-time (expressed in seconds) in relation to the size $N \times M$, and to $w>0$, for the two approaches. 
  
  From the results of Table \ref{tab_times_1}, Table \ref{tab_times_2}, and Table \ref{tab_times}, it is clear that the approach (2) is faster than (1). 
\begin{table}[!h]
\begin{center}
\begin{tabular}{ccccccccccc}
\hline
{\bf N x M} & \textbf{w=1} & \textbf{w=4} & \textbf{w=9} & \textbf{w=25} & \textbf{w=100} & \textbf{w=400}\\
\hline
$1 \times 1$	& 0.041870	& 0.040219	& 0.041105	& 0.038928	& 0.042712	& 0.039214\\
$2 \times 2$	&0.042075	&0.042040	&0.041391	&0.041092	&0.040991	&0.042992\\
$3 \times 3$	&0.043006	&0.042132	&0.041864	&0.044914	&0.043080	&0.048435\\
$5 \times 5$	&0.044211	&0.048586	&0.047370	&0.048077	&0.055869	&0.060404\\
$10 \times 10$	&0.064676	&0.066046	&0.068574	&0.074969	&0.093118	&0.227232\\
\hline 
\end{tabular}
\caption{{\small CPU-time for the data sets of the dimension listed in the first column, corresponding to each value of $w$, for the approach (1).}} \label{tab_times_1}                          
\end{center}
\end{table}
\begin{table}[!h]
\begin{center}
\begin{tabular}{ccccccccccc}
\hline
 {\bf N x M} &  \textbf{w=1} & \textbf{w=4} & \textbf{w=9} & \textbf{w=25} & \textbf{w=100} & \textbf{w=400}\\
\hline
$1 \times 1$ &0.043561&	0.040653&	0.042252&	0.043928&	0.040932	&0.041430\\
$2 \times 2$ &0.043919&	0.044333&	0.043065&	0.043307&	0.043592	&0.045282\\
$3 \times 3$ &0.043289&	0.043140&	0.044113&	0.043985&	0.041659	&0.044925\\
$5 \times 5$ &0.043000&	0.043482&	0.043960&	0.044029&	0.044554	&0.045373\\
$10 \times 10$ &0.045015&	0.045146&	0.045300&	0.046053&	0.058119	&0.147052\\
\hline 
\end{tabular}
\caption{{\small CPU-time for the data sets of the dimension listed in the first column, corresponding to each value of $w$, for the approach (2).}} \label{tab_times_2}                          
\end{center}
\end{table}
\begin{figure}
\centering
\includegraphics[scale=0.7]{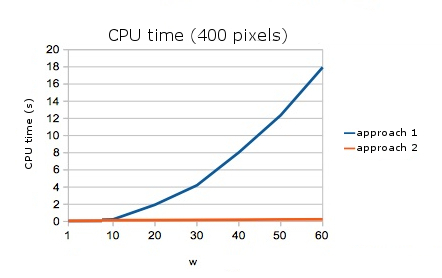}
\caption{Comparison between the CPU-time for the approaches (1) and (2), in function of $w>0$ for an image of size $20 \times 20$.}\label{fig10}
\end{figure}

The advantage provided by procedure (2) appears more evident if we consider matrices of higher dimension with respect to those considered in Table \ref{tab_times_1} and Table \ref{tab_times_2}. Indeed, in Table \ref{tab_times} we tested a matrix of dimension $20 \times 20$, with various values of $w$.
\begin{table}[!h]
\begin{center}
\begin{tabular}{ccccccccccc}
\hline
{\bf Approach} & \textbf{w=10} & \textbf{w=20} & \textbf{w=30} & \textbf{w=40} & \textbf{w=50} & \textbf{w=60} \\
\hline
(1) &	0.236207 &	1.943428 &	4.20498 & 	8.0352 &	12.3671 &	17.9449 \\
(2) &	0.14637 & 	0.158876 &	0.18029 &	0.20343 &	0.23314 &	0.24744 \\
\hline 
\end{tabular}
\caption{{\small CPU-time for a matrix of size $20 \times 20$ for both approaches (1) and (2), for several values of $w$.}}   \label{tab_times}                      
\end{center}
\end{table}   

  The results of Table \ref{tab_times} have been plotted in Fig. \ref{fig10}; the execution times using the algorithm based on (1) is longer than the corresponding ones based on (2).
           
  On the basis of the previous considerations, approach (2) results more efficient. For these reasons, the choice to implement the second procedure has been preferred.


\section{Discussion of results and conclusions} \label{sec5} 

\begin{figure}
\centering
\includegraphics[scale=0.8]{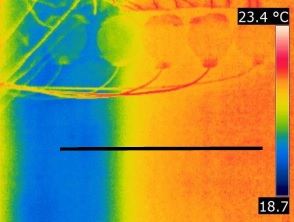}
\hskip0.6cm
\includegraphics[scale=0.8]{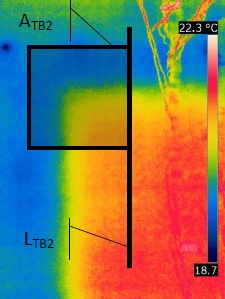}
\caption{The black line has been used for the computation of the index $I_{tb}$ in the 2D thermal bridge.} \label{fig12}
\end{figure}

The segmentation method developed in Section \ref{sec2} and Section \ref{sec3} is applied in order to detect the shape of thermal bridges of the building envelope from thermographic images. Moreover, it also allows to determine their heat losses, using a suitable {\em incidence factor of thermal bridge}, i.e., the index $I_{tb}$, previously introduced in \cite{ASBA1,ASBA2}:
\be
I_{tb}\ :=\ \frac{\displaystyle \sum_{p=1}^N \, (T_i - T_{p})}{\displaystyle N (T_i - T_{1D})},
\ee
where $T_i$ is the internal air temperature, $T_{p}$ is the acquired surface temperature of the single pixel from the infrared camera, $T_{1D}$ is the surface temperature of the undisturbed zone of the wall, evaluated with infrared camera as well, and $N$ is the number of the pixels that compose an imaginary line (see e.g., Fig. \ref{fig12}) along the thermal bridge. The $I_{tb}$ is a parameter higher than one and it gives information on the thermal bridge effects on the overall energy performance of the investigated wall. The accuracy of $I_{tb}$ depends on the resolution of the thermal map, and obviously it is independent by the length of the imaginary line chosen for its computation.

   The index $I_{tb}$ has been calculated using the original thermographic images and by the corresponding image reconstructed by the application of the S-K algorithm. In the first case, the shape of the thermal bridge has been manually identified, while in the second case it has been automatically extracted by the proposed segmentation algorithm. The above achieved results have been compared with a reference index computed with the real data coming from the experimental setup, i.e., with the thermal parameters detected by heat flow meters and thermocouples in the hot-box setup with controlled laboratory conditions. We can state that the $I_{tb}$ index computed thanks to the application of the segmentation algorithm developed in this paper provides results that are closer to the reference $I_{tb}$ index derived from sensor measurements that directly acquire the thermal flux. 
     
   More precisely, in case of the examined typologies of thermal bridges, the values of the index $I_{tb}$ determined by the proposed procedure and that one obtained by the experimental sensor methodology, show a significant improvement.
Here, the $I_{tb}$ indexes computed from the original method described in \cite{ASBA1}, i.e., from the original not-enhanced thermographic images, are $1.611$ and $1.467$, respectively. While, the corresponding ones computed thanks to the application of the proposed algorithm, are $1.585$ and $1.462$, respectively. Thus, the latter are closer to the reference values of $I_{tb}$ (obtained by probes), which are $1.439$ and $1.303$, respectively. Indeed, we point out that if we estimate the ratio between the absolute errors of the two different methods compared with the reference value of $I_{tb}$, it turns out that the new proposed approach achieves an improvement of around $15 \%$ and $4\%$, respectively. Obviously, this reflects in a more accurate heat loss estimation. In the above examples, one of the main aspects that allows to achieve such improvement is given by the fact that the thermal bridges are generated by different materials, then, the proposed segmentation algorithm detects with a high accuracy the boundary which separates the two different areas. 

  In conclusion, the above experimental results show the performance and the enhancement of the procedure introduced in this paper with respect to the standard approach for the evaluation of the heat losses in buildings. In addition, further advantages of the present algorithm consist on the possibility to detect the exact geometry of the thermal bridge, and to be automatic (not operator depending).


\vskip0.2cm

\section*{Acknowledgments}

The work has been supported by the research project: ``Enhancement algorithms of thermographic images to study the influence of thermal bridges in the energetic analysis of buildings'' funded by the ``Fondazione Cassa di Risparmio di Perugia''. The authors D. Costarelli, M. Seracini, and G. Vinti, are members of the Dipartimento di Matematica e Informatica of the University of Perugia and of the GNAMPA of INdAM (Istituto Nazionale di Alta Matematica).  Moreover, D. Costarelli holds a research grant (Post-Doc) funded by the INdAM, and together with M. Seracini, they are partially supported by the 2017 GNAMPA-INdAM project: ``Approssimazione con operatori discreti e problemi di minimo per funzionali del calcolo delle variazioni con applicazioni all'imaging''.

\vskip0.2cm

%


\begin{thebibliography}{99}

\bibitem{AGBA} P.N.~Agrawal, B.~Baxhaku, Degree of approximation for bivariate extension of Chlodowsky-type q-Bernstein-Stancu-Kantorovich operators. {\em Applied Mathematics and Computation}, 306, 2017, pp. 56--72.

\bibitem{ANVI1} L.~Angeloni and G.~Vinti, Approximation in variation by homothetic operators in multidimensional setting. {\em Differential and Integral Equations} 26 (5-6) 2013, pp. 655--674.

\bibitem{ASBA1} F. Asdrubali, G. Baldinelli, Thermal transmittance measurements with the hot box method: calibration, experimental procedures, and uncertainty analyses of three different approaches. {\em Energy Build} 43, 2011, pp. 1618--26. 

\bibitem{ASBA2}  F. Asdrubali, G. Baldinelli, F. Bianchi, A quantitative methodology to evaluate thermal bridges in buildings. {\em Applied Energy} 97, 2012, pp. 365--373.

\bibitem{ABBP1} F. Asdrubali, G. Baldinelli, F. Bianchi, A.L. Pisello, Infrared thermography assessment of thermal bridges in building envelope: Experimental validation in a test room setup. {\em Sustainability} 6(10), 2014, pp. 7107--7120.

\bibitem{BBPR1} G. Baldinelli, F. Bianchi, A. Presciutti, A. Rotili, Transient Heat Transfer in Radiant Floors: A Comparative Analysis between the Lumped Capacitance Method and Infrared Thermography Measurements. {\em J. Imaging} 2(3) 22, 2016.

\bibitem{BABUSTVI2} C.~Bardaro, P.L.~Butzer, R.L.~Stens and G.~Vinti, Kantorovich-type generalized sampling series in the setting of Orlicz spaces. {\em Sampl.~Theory~Signal Image~Process.} 6 (1), 2007, pp.~29--52.
 
\bibitem{BABUSTVI3} C.~Bardaro, P.L.~Butzer, R.L.~Stens and G.~Vinti, Prediction by samples from the past with error estimates covering discontinuous signals. {\em IEEE~Trans.~Inform.~Theory} 56 (1), 2010, pp.~614--633.

\bibitem{BAMUVI} C.~Bardaro, J.~Musielak and G.~Vinti, {\em Nonlinear Integral Operators and Applications,} De Gruyter Series in Nonlinear Analysis and Applications 9, New York--Berlin 2003.

\bibitem{BOCAME1} S. Boccardi, G.M. Carlomagno, C. Meola, An Excursus on Infrared Thermography Imaging, {\em Journal of Imaging} DOI: 10.3390/jimaging2040036, 2016. 

\bibitem{BOCAME2} S. Boccardi, G.M. Carlomagno, C. Meola, The added value of infrared thermography in the measurement of temperature-stress coupled effects, {\em Sensors and Transducers}, 2016.

\bibitem{BUFIST} P.L.~Butzer, A.~Fisher and R.L.~Stens, Generalized sampling approximation of multivariate signals: theory and applications. {\em Note di Matematica}, 10 (1), 1990, pp.~173--191.

\bibitem{BUNE} P.L.~Butzer and R.J.~Nessel, {\em Fourier Analysis and Approximation I,} Academic Press, New York--London 1971.

\bibitem{BURIST2} P.L.~Butzer, S.~Ries and R.L.~Stens,  
Approximation of continuous and discontinuous functions by generalized sampling series. {\em J.~Approx.~Theory} 50, 1987, pp.~25--39.

\bibitem{CACA1} G.M. Carlomagno, G. Cardone, Infrared thermography for convective heat transfer measurements. {\em Experiments in Fluids}, 49 (6), 2010, 1187-1218.

\bibitem{CHYUVE1} S.G. Chang, B. Yu, B., and M. Vetterli, Adaptive wavelet thresholding for image denoising and compression. {\em IEEE trans. image proc.}, 9(9), 2000, pp.~1532--1546.

\bibitem{CHDI1} C.K. Chui, and H. Diamond, A natural formulation of quasi-interpolation by multivariate splines. {\em Proc. American Math. Society} 99(4), 1987, pp.~643--646.

\bibitem{CLCOMIVI1} F.~Cluni, D.~Costarelli, A.M.~Minotti and G.~Vinti, Applications of sampling Kantorovich operators to thermographic images for seismic engineering, {\em J.~Comput.~Anal.~Appl.} 19 (4), 2015, pp.~602--617.

\bibitem{CLCOMIVI2} F.~Cluni, D.~Costarelli, A.M.~Minotti and G.~Vinti, Enhancement of thermographic images as tool for structural analysis in earthquake engineering. {\em NDT \& E International} 70, 2015, pp.~60--72. 

\bibitem{COGA1} L. Coroianu, and S.G. Gal, Approximation by max-product sampling operators based on sinc-type kernels. {\em Sampling Th. Signal Image Processing} 10 (3), 2011, pp.~211-230.

\bibitem{COGA2} L. Coroianu, and S.G. Gal, Saturation results for the truncated max-product sampling operators based on sinc and Fej\'er-type kernels. {\em Sampling Th. Signal Image Processing} 11 (1), 2012, pp.~113-132. 

\bibitem{CO2} D.~Costarelli, Neural network operators: interpolation of multivariate functions. {\em Neural Networks} 67, 2015, pp. 28--36.

\bibitem{COMIVI1} D. Costarelli, A.M. Minotti, G. Vinti, Approximation of discontinuous signals by sampling Kantorovich series. {\em J. Math. Anal. Appl.} 450 (2), 2017, pp. 1083--1103.

\bibitem{COVI1} D.~Costarelli and G.~Vinti, Approximation by multivariate generalized sampling Kantorovich operators in the setting of Orlicz spaces, {\it Bollettino~U.M.I.} 4 (9), 2011, pp.~445--468.

\bibitem{COVI2} D.~Costarelli and G.~Vinti, Approximation by nonlinear multivariate sampling Kantorovich type operators and applications to image processing, {\em Num.~Funct.~Anal.~Opt.} 34 (8), 2013, pp.~819--844.

\bibitem{COVI4} D.~Costarelli and G.~Vinti, Rate of approximation for multivariate sampling Kantorovich operators on some functions spaces, {\em J.~Integral~Eq.~Appl.}, 26 (4), 2014, pp. 455--481.

\bibitem{COVI5} D.~Costarelli and G.~Vinti, Degree of approximation for nonlinear multivariate sampling Kantorovich operators on some functions spaces. {\em Num. Funct. Anal. Opt.} 36 (8), 2015, pp. 964--990. 

\bibitem{COVI6} D.~Costarelli and G.~Vinti, Approximation by max-product neural network operators of Kantorovich type. {\em Results in Mathematics} 69 (3), 2016, pp. 505--519. 

\bibitem{COVI7} D.~Costarelli and G.~Vinti, Max-product neural network and quasi-interpolation operators activated by sigmoidal functions. {\em J. Approx. Theory} 209, 2016, pp. 1--22.

\bibitem{COVI8} D.~Costarelli and G.~Vinti, Pointwise and uniform approximation by multivariate neural network operators of the max-product type. {\em Neural Networks} 81, 2016, pp. 81--90.

\bibitem{COVI9} D. Costarelli, G. Vinti, Convergence for a family of neural network operators in Orlicz spaces. {\em Mathematische Nachrichten} 290 (2-3), 2017, pp. 226--235.  

\bibitem{GOWO} R.~Gonzales and R.~Woods, {\em Digital Image Processing}, Prentice-Hall NJ - USA, 2002.

\bibitem{IC} G. I\c{c}oz, A Kantorovich variant of a new type Bernstein-Stancu polynomials. {\em Applied Mathematics and Computation} 218(17), 2012, pp. 8552--8560.

\bibitem{KASAWO1} J.N Kapur, P.K. Sahoo, and A.K. Wong, A new method for gray-level picture thresholding using the entropy of the histogram. {\em Computer vision, graphics, and image processing}, 29(3), 1985, pp.~273--285.

\bibitem{KI} M.~Kass, A.~Vitkin and D.~Terzopoulos, Snakes: Active Contour Models, {\em WSEAS~Trans.~on Signal~Processing}, 10, 2014, pp.~288--300.

\bibitem{KIME1} A. Kivinukk, T. Metsmagi, On boundedness inequalities of some semi-discrete operators in connection with sampling operators. {\em 2015 International Conference on Sampling Theory and Applications, SampTA 2015}, 2015, pp.~48--52.

\bibitem{KITA1} A. Kivinukk, and G. Tamberg, Interpolating generalized Shannon sampling operators, their norms and approximation properties. {\em Sampling Theory in Signal and Image Processing}, 8, 2009, pp.~77--95.

\bibitem{KLME} E. Klein, H. J. Metz, P. Stucki, {\em Advances in Digital Image Processing  Theory, Application, Implementation}, Springer - USA, 1979.

\bibitem{KV1}  B.I. Kvasov, Methods of shape-preserving spline approximation. {\em Singapore: World Scientific} 2000.

\bibitem{LIWA1} C.J. Li, and R.H. Wang, Bivariate quartic spline spaces and quasi-interpolation operators. {\em J. Comput. Applied Math.}, 190(1), 2006, pp.~325--338.

\bibitem{MU1} J.~Musielak, {\it Orlicz Spaces and Modular Spaces,} Springer-Verlag, Lecture~Notes~Math., 1983.

\bibitem{ORTA1} O. Orlova, and G. Tamberg, On approximation properties of Kantorovich-type sampling operators. {\em Sam. Th. Signal Image Proc.} 201, 2016, pp.~73--86.

\bibitem{SA1} D. Salomon, B-Spline Approximation. {\em The Computer Graphics Manual. Springer London}, 2011. 

\bibitem{SK} V.~Skala, Fast interpolation and approximation of scattered multidimensional and dynamic data using radial basis functions, {\it WSEAS~Trans.~on~Mathematics}, 12 (5), 2013, pp.~501--511.

\bibitem{TA1} G. Tamberg, On truncation errors of some generalized Shannon sampling operators. {\em Numerical Algorithms}, 55(2), 2010, pp.~367--382.

\bibitem{VE1} R. Verfurth, Error estimates for some quasi-interpolation operators. {\em ESAIM: Mathematical Modelling and Numerical Analysis} 33(4), 1999, pp.~695--713.

\bibitem{VIZA1} G.~Vinti and L.~Zampogni, Approximation by means of nonlinear Kantorovich sampling type operators in Orlicz spaces, {\em J.~Approx.~Theory} 161, 2009, pp.~511--528.

\bibitem{VIZA3} G.~Vinti, and L.~Zampogni, Approximation results for a general class of
Kantorovich type operators, {\em Adv.~Nonlin.~Studies} 14, 2014, pp.~991--1011.

\bibitem{ZH1}  Y. J. Zhang, A survey on evaluation methods for image segmentation. {\em Pattern recognition}, 29(8), 1996, pp.~1335--1346.

\end{thebibliography}
\end{document}